\documentclass{elsarticle2}

\usepackage{amsmath}
\usepackage{amsthm}
\usepackage{amsfonts}
\usepackage{tikz}
\usepackage{enumerate}
\usepackage{hyperref}
\usetikzlibrary{patterns}

\title{A Formula for the M\"obius function of the Permutation Poset Based on a Topological Decomposition}
\author{Jason P Smith}
\date{}
\address{Department of Computer and Information Sciences\\ University of Strathclyde, Glasgow, UK}

\newtheorem{thm}{Theorem}
\newtheorem{lem}[thm]{Lemma}
\newtheorem{ex}[thm]{Example}

\newtheorem{conj}[thm]{Conjecture}
\newtheorem{prop}[thm]{Proposition}
\newtheorem{defn}[thm]{Definition}
\newtheorem{rem}[thm]{Remark}

\newtheorem{que}[thm]{Question}

\DeclareMathOperator{\NE}{\text{\rm NE}}
\DeclareMathOperator{\EZ}{\text{\rm EZ}}

\begin{document}
\begin{abstract}
We present a two term formula for the M\"obius function of intervals in the poset of all permutations, ordered by pattern containment.  The first term in this formula is the number of so called normal occurrences of one permutation in another.  Our definition of normal occurrences is similar to those that have appeared in several variations in the literature on the M\"obius function of this and other posets, but simpler than most of them.  The second term in the formula is complicated, but we conjecture that it equals zero for a significant proportion of intervals. We present some cases where the second term vanishes and others where it is nonzero. Computing the M\"obius function recursively from its definition has exponential complexity, whereas the computation of the first term in our formula is polynomial and the exponential part is isolated to the second term, which seems to often vanish.  We also present a result on the M\"obius function of posets connected by a poset fibration.
\end{abstract}

	\maketitle

\section{Introduction}
Let $\sigma$ and $\pi$ be permutations of positive integers. We define an \emph{occurrence} of $\sigma$ in $\pi$ to be a subsequence of $\pi$ with the same relative order of size as the letters in $\sigma$. For example, $132$ occurs twice in $23541$, as the subsequences 254 and 354. The \emph{permutation poset} $\mathcal{P}$ consists of all permutations with the partial order $\sigma\le\pi$ if there is an occurrence of $\sigma$ in $\pi$. An \emph{interval} $[\sigma,\pi]$ in $\mathcal{P}$ is the subposet $\{z\in\mathcal{P}\,|\,\sigma\le z\le\pi\}$. The \emph{M\"obius function} for a poset is defined recursively as: $\mu(a,b)=0$ if $a\not\le b$, $\mu(a,a)=1$ for all $a$ and, for $a<b$: 
$$\mu(a,b)=-\sum_{a\le z< b}\mu(a,z).$$ 

The first systematic study of the M\"obius function of general posets appeared in \cite{Rota64} and the first result pertaining to the M\"obius function of intervals of $\mathcal{P}$ appeared in \cite{SagVat06}, where a formula for intervals of  \emph{layered} permutations was presented. A layered permutation is the direct sum of decreasing permutations, where the \emph{direct sum} $\sigma\oplus\pi$ of two permutations $\sigma$ and~$\pi$ is obtained by appending~$\pi$ to $\sigma$ after adding the length of $\sigma$ to each letter of $\pi$. For example,~$312\oplus213=312546$. There is an analogous \emph{skew sum} $\sigma\ominus\pi$ where $\pi$ is appended to $\sigma$ after the length of $\pi$ is added to each element of $\sigma$. In~\cite{BJJS11} a formula for the M\"obius function is presented for intervals of  \emph{decomposable} permutations, that is, permutations that can be written as the direct sum of two or more non-empty permutations. This formula, however, is recursive and bottoms out in intervals bounded by indecomposable permutations, for which there is no general formula for the M\"obius function.

Furthermore, in \cite{BJJS11} a formula is presented for intervals of  \emph{separable} permutations, that is, permutations that avoid 2413 and 3142, or equivalently, permutations that can be written using only direct sums, skew sums and the singleton permutation 1. A formula for the M\"obius function of intervals of permutations with a fixed number of descents is given in \cite{Smith14}, where a \emph{descent} occurs at position~$i$ in a permutation~$\pi=\pi_1\ldots\pi_n$  if $\pi_i>\pi_{i+1}$. Further results have been presented in~\cite{McSt13,Smith13,SteTen10}. However, the proportion of intervals $[\sigma,\pi]$ which satisfy any of these properties approaches zero as the length of~$\pi$ increases. There are indications that the formula we present here reduces the computation of the M\"obius function to polynomial time for a significant proportion of intervals.

Many of the results on the M\"obius function of intervals of $\mathcal{P}$, and also of some posets of words, are linked to the number of what have been termed \emph{normal occurrences}, or \emph{normal embeddings}, in the literature, see \cite{Bjo90,Bjo93,BJJS11,SagVat06,Smith14}. The first appearance  of normal occurrences is in Bj\"orner's paper \cite{Bjo90} where a formula for the M\"obius function of intervals of words with subword order is presented. The definition of a normal occurrence has varied in these papers, but all follow a similar theme. 

Our definition of normal occurrences, which is simpler than most previous ones, is based upon the \emph{adjacencies} of a permutation, where an adjacency in a permutation is a maximal sequence of increasing or decreasing consecutively valued letters in consecutive positions and the \emph{tail} of an adjacency is all but its first letter. A \emph{normal occurrence} of $\sigma$ in $\pi$, in our definition, is any occurrence that includes all the tails of all the adjacencies of $\pi$. This definition of normal occurrences based on adjacencies does not seem to have been considered previously, but in~\cite{Smith14} we presented a slightly different version.

We present a formula, in Theorem \ref{thm:main}, that shows the M\"obius function of~$[\sigma,\pi]$ is, up to a sign, equal to the number of normal occurrences of $\sigma$ in~$\pi$ plus an extra term that seems to vanish for a significant proportion of intervals. For example, we know this extra term vanishes if $\sigma$ and $\pi$ have the same number of descents, which is a consequence of the result in \cite{Smith14}. Using \emph{interval blocks}, which appear in  \cite{SteTen10}, we prove that if for all permutations~$\lambda\in[\sigma,\pi)$ there is a singleton interval block, that is, a letter of $\pi$ which belongs to no occurrence of $\lambda$, the second term of the formula vanishes. The above mentioned cases are of zero proportion when the length of $\pi$ goes to infinity, but computer tests indicate that for a substantial proportion of intervals the second term of our formula vanishes. Why that is the case is still a mystery, but this suggests that many more families of intervals than are now known may turn out to have a tractable M\"obius function.

It is shown in \cite{McSt13} that if $\pi$ is decomposable and has equal consecutive components then for any subpermutation $\sigma$ obtained by removing $k>1$ of the equal components, the interval $[\sigma,\pi]$ contains a disconnected subinterval. Many of the definitions of normal occurrences have an extra condition for the case when~$\pi$ has this property. We prove a result that indicates the second term of our formula for the M\"obius function is often non-zero for the slightly simpler case that $\pi$ is decomposable and all the components are equal. Exactly what the connection is between this second term and the topology of such intervals is another mystery.

Computing the M\"obius function using the original recursive formula has exponential complexity, whereas our formula splits the computation into two parts. The first part, that is, computing the number of normal occurrences, can be done in polynomial time and the second part has exponential complexity in the general case, but computational evidence suggests that in a significant proportion of cases this second term vanishes. Our formula here is the first formula for arbitrary intervals of permutations that seems to have polynomial time complexity for a significant proportion of intervals. 

In Section \ref{sec:defn} we introduce some definitions, give a brief introduction to the topology of posets and present a poset fibration of $[\sigma,\pi]$ that we later use to compute $\mu(\sigma,\pi)$. In Section \ref{sec:result} we present and prove our main result, that the M\"obius function of intervals of $\mathcal{P}$ equals the number of normal occurrences plus an extra term that we define. In Section \ref{sec:fibration} we present a result that links the M\"obius function of two posets connected by a poset fibration satisfying a certain condition. This indicates there is possibly a more general condition for the main result of Bj\"orner, Wachs and Welker in~\cite{Bjo05}. In Section \ref{sec:app} we apply our formula to show that the  M\"obius function of~$[\sigma,\pi]$ equals the number of normal occurrences of $\sigma$ in~$\pi$ if for each $\lambda\in[\sigma,\pi)$ there is at least one letter of $\pi$ which is not in any occurrence of $\lambda$. We also show that the value of the second term of our formula for the M\"obius function of $[\sigma,\pi]$ is often nonzero when~$\pi$ can be decomposed into the direct sum of equal components. Furthermore, we consider for which permutations all occurrences are normal.

\section{Definitions and Preliminaries}\label{sec:defn}
In this section we introduce some definitions required to present our main result. We begin with an important property of permutations that is fundamental to our results:
\begin{defn}\label{defn:emb}
An \emph{adjacency} in a permutation  is a maximal sequence, of \linebreak length $\ell\ge1$, of increasing or decreasing consecutively valued letters in consecutive order. The \emph{tail} of an adjacency of length at least 2 is all but the first letter of the adjacency. An adjacency of length 1 does not have a tail.
\end{defn}

\begin{ex}
The permutation $\pi=2314765$ has adjacencies 23, 1, 4 and 765 and the tails are 3 and~65.
\end{ex}

Next we define embeddings and our version of normal embeddings. Embeddings are in one-to-one correspondence with occurrences, and we use embeddings instead of occurrences throughout the rest of the paper because they allow for easier presentation of the required definitions. 
\begin{defn}\label{defn:normal}
Consider permutations $\sigma\le\pi$.  An \emph{embedding} $\eta$ of $\sigma$ in $\pi$ is a sequence of the same length as $\pi$ such that the nonzero letters in $\eta$ are the letters of an occurrence of $\sigma$ in $\pi$ and in the same positions in $\eta$ as in $\pi$.

An embedding $\eta$ of $\sigma$ in $\pi$ is \emph{normal} if the positions of all the letters in all the tails of the adjacencies in $\pi$ are nonzero in $\eta$. We denote the number of normal embeddings of $\sigma$ in $\pi$ as~$\NE(\sigma,\pi)$.
\end{defn}

\begin{ex}
For $\sigma=132$ and $\pi=2314765$ the sequence $0300065$ is the only normal embedding of $\sigma$ in $\pi$, so $\NE(\sigma,\pi)=1$.
\end{ex}

\begin{prop}\label{prop:complexity}
Computing $\NE(\sigma,\pi)$ for a fixed $\sigma$ can be done in time polynomial in the length of $\pi$.
\begin{proof}
Counting the number of occurrences of $\sigma$ in $\pi$, of lengths $k$ and $n$, respectively, can be done in polynomial time $\mathcal{O}(n^k)$ by exhaustive search, and testing for normality is linear.
\end{proof} 
\end{prop}

We use the adjacencies of a permutation to break down the permutation and embeddings into smaller components.
\begin{defn}
Consider permutations $\sigma\le\pi$ and an embedding $\eta$ of $\sigma$ in $\pi$. Let $\hat{\pi}=(\hat{\pi}_1,\ldots,\hat{\pi}_t)$ be the decomposition of $\pi$ into its adjacencies, that is,~$\hat{\pi}_i$ is a maximal increasing or decreasing permutation corresponding to the $i$'th adjacency of $\pi$. 

Define $\hat{\eta}:=(\hat{\eta}_1,\ldots,\hat{\eta}_t)$ where $\hat{\eta}_i$ is the permutation obtained from the non\-zero letters that $\eta$ embeds in the $i$'th adjacency of $\pi$. If $\eta$ does not embed in any letters of the $i$'th adjacency then $\hat{\eta}_i=\emptyset$.
\end{defn}

\begin{ex}
If $\sigma=132$ and $\pi=2314765$ then $\hat{\pi}=(12,1,1,321)$ and the embedding $\eta=0010760$ gives $\hat{\eta}=(\emptyset,1,\emptyset,21)$.
\end{ex}

When considering  embeddings  the selection of letters within an adjacency is usually irrelevant. This is made formal by the following equivalence relation.
\begin{defn}
Let $E^{\sigma,\pi}$ be the set of embeddings of $\sigma$ in $\pi$. Define an equivalence relation on embeddings where $\eta\sim\psi$ if the only differences between $\eta$ and $\psi$ occur within adjacencies of $\pi$. Define $\widehat{E}^{\sigma,\pi}$ as the set containing the rightmost embedding, that is, the embedding where the nonzero letters are the furthest right, of each equivalence class of $E^{\sigma,\pi}\slash\sim$.

Consider $\eta\in \widehat{E}^{\sigma,\pi}$ and define the zero set of $\eta$ as $Z(\eta)=\{i\,|\, \eta_i=0\}$. Define $\EZ^{\sigma,\pi}$ to be the set of sets of embeddings in $\widehat{E}^{\sigma,\pi}$ such that for each set~$S\in \EZ^{\sigma,\pi}$ we have $\bigcap_{\eta\in S}Z(\eta)=\emptyset$.
\end{defn}

When defining $\widehat{E}^{\sigma,\pi}$ we choose the rightmost embedding to ensure that all normal embeddings are in $\widehat{E}^{\sigma,\pi}$. Note that if $\eta\sim\psi$ then $\hat{\eta}=\hat{\psi}$, which can be used as an equivalent definition of the equivalence relation. The set $\EZ^{\sigma,\pi}$ is upwards closed under containment because if we take any set $S\in \EZ^{\sigma,\pi}$ adding a new embedding to $S$ will result in a set that still has empty intersection of zero sets. 

\begin{ex}\label{ex:3}
If $\sigma=132$ and $\pi=413265$ then the embedding $013200$ has zero set $Z(013200)=\{1,5,6\}$ and \begin{align*}E^{\sigma,\pi}=&\{013200,\,400065,\,010065,\,003065,\,000265\},\\ \widehat{E}^{\sigma,\pi}=&\{013200,\,400065,\,010065,\,000265\}, \\\EZ^{\sigma,\pi}=&\{\{013200,\,400065\},\{013200,\,400065,\,010065\},\\ &\{013200,\,400065,\,000265\},\{013200,\,400065,\,010065,\,000265\}\}.\end{align*}
\end{ex}

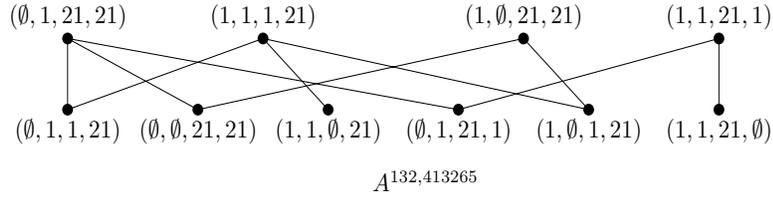
\begin{figure}\centering\resizebox{345pt}{400pt}{\begin{tikzpicture}
\draw (0,0) node[fill,circle,scale=0.5pt]{};
\draw (0,0) node[below] {$(\emptyset,1,21,\emptyset)$};
\draw (-1,1) node[fill,circle,scale=0.5pt]{};
\draw (-1,1) node[left] {$(1,1,21,\emptyset)$};
\draw (1,1) node[fill,circle,scale=0.5pt]{};
\draw (1,1) node[right] {$(\emptyset,1,21,1)$};
\draw (-1,2) node[fill,circle,scale=0.5pt]{};
\draw (-1,2) node[left] {$(1,1,21,1)$};
\draw (1,2) node[fill,circle,scale=0.5pt]{};
\draw (1,2) node[right] {$(\emptyset,1,21,21)$};
\draw (0,3) node[fill,circle,scale=0.5pt]{};
\draw (0,3) node[above] {$(1,1,21,21)$};
\draw (0,0) -- (-1,1) -- (-1,2) -- (0,3) -- (1,2) -- (1,1) -- (0,0);
\draw (-1,2) -- (1,1);
\draw (0,-1) node{$P(013200)=[\emptyset,1]\times[1,1]\times[21,21]\times[\emptyset,21]$};
\draw (7,0) node[fill,circle,scale=0.5pt]{};
\draw (7,0) node[below] {$(1,\emptyset,\emptyset,21)$};
\draw (6,1) node[fill,circle,scale=0.5pt]{};
\draw (6,1) node[left] {$(1,1,\emptyset,21)$};
\draw (8,1) node[fill,circle,scale=0.5pt]{};
\draw (8,1) node[right] {$(1,\emptyset,1,21)$};
\draw (6,2) node[fill,circle,scale=0.5pt]{};
\draw (6,2) node[left] {$(1,1,1,21)$};
\draw (8,2) node[fill,circle,scale=0.5pt]{};
\draw (8,2) node[right] {$(1,\emptyset,21,21)$};
\draw (7,3) node[fill,circle,scale=0.5pt]{};
\draw (7,3) node[above] {$(1,1,21,21)$};
\draw (7,0) -- (6,1) -- (6,2) -- (7,3) -- (8,2) -- (8,1) -- (7,0);
\draw (6,2) -- (8,1);
\draw (7,-1) node {$P(400065)=[1,1]\times[\emptyset,1]\times[\emptyset,21]\times[21,21]$};
\draw (0,-6) node[fill,circle,scale=0.5pt]{};
\draw (0,-6) node[below] {$(\emptyset,1,\emptyset,21)$};
\draw (-1,-5) node[fill,circle,scale=0.5pt]{};
\draw (-1,-5) node[left] {$(1,1,\emptyset,21)$};
\draw (1,-5) node[fill,circle,scale=0.5pt]{};
\draw (1,-5) node[right] {$(\emptyset,1,1,21)$};
\draw (-1,-4) node[fill,circle,scale=0.5pt]{};
\draw (-1,-4) node[left] {$(1,1,1,21)$};
\draw (1,-4) node[fill,circle,scale=0.5pt]{};
\draw (1,-4) node[right] {$(\emptyset,1,21,21)$};
\draw (0,-3) node[fill,circle,scale=0.5pt]{};
\draw (0,-3) node[above] {$(1,1,21,21)$};
\draw (0,-6) -- (-1,-5) -- (-1,-4) -- (0,-3) -- (1,-4) -- (1,-5) -- (0,-6);
\draw (-1,-4) -- (1,-5);
\draw (0,-7) node {$P(010065)=[\emptyset,1]\times[1,1]\times[\emptyset,21]\times[21,21]$};
\draw (7,-6) node[fill,circle,scale=0.5pt]{};
\draw (7,-6) node[below] {$(\emptyset,\emptyset,1,21)$};
\draw (5.5,-5) node[fill,circle,scale=0.5pt]{};
\draw (5.5,-5) node[left] {$(1,\emptyset,1,21)$};
\draw (8.5,-5) node[fill,circle,scale=0.5pt]{};
\draw (8.5,-5) node[right] {$(\emptyset,1,1,21)$};
\draw (7,-5) node[fill,circle,scale=0.5pt]{};
\draw (7,-5) node[below] {$(\emptyset,\emptyset,21,21)$};
\draw (7,-4) node[fill,circle,scale=0.5pt]{};
\draw (7,-4) node[above] {$(1,1,1,21)$};
\draw (5.5,-4) node[fill,circle,scale=0.5pt]{};
\draw (5.5,-4) node[left] {$(1,\emptyset,21,21)$};
\draw (8.5,-4) node[fill,circle,scale=0.5pt]{};
\draw (8.5,-4) node[right] {$(\emptyset,1,21,21)$};
\draw (7,-3) node[fill,circle,scale=0.5pt]{};
\draw (7,-3) node[above] {$(1,1,21,21)$};
\draw (7,-6) -- (5.5,-5) -- (7,-4) -- (7,-3) -- (8.5,-4) -- (8.5,-5) -- (7,-6);
\draw (7,-6) -- (7,-5);
\draw (7,-5) -- (5.5,-4) -- (7,-3);
\draw (5.5,-5) -- (5.5,-4);
\draw (7,-5) -- (8.5,-4);
\draw (7,-4) -- (8.5,-5);
\draw (7,-7) node {$P(000265)=[\emptyset,1]\times[\emptyset,1]\times[1,21]\times[21,21]$};
\draw (-2,-9) node[fill,circle,scale=0.5pt]{};
\draw (-2,-9) node[above] {$(\emptyset,1,21,21)$};
\draw (1,-9) node[fill,circle,scale=0.5pt]{};
\draw (1,-9) node[above] {$(1,1,1,21)$};
\draw (5,-9) node[fill,circle,scale=0.5pt]{};
\draw (5,-9) node[above] {$(1,\emptyset,21,21)$};
\draw (8,-9) node[fill,circle,scale=0.5pt]{};
\draw (8,-9) node[above] {$(1,1,21,1)$};
\draw (-2,-10) node[fill,circle,scale=0.5pt]{};
\draw (-2,-10) node[below] {$(\emptyset,1,1,21)$};
\draw (-0,-10) node[fill,circle,scale=0.5pt]{};
\draw (-0,-10) node[below] {$(\emptyset,\emptyset,21,21)$};
\draw (2,-10) node[fill,circle,scale=0.5pt]{};
\draw (2,-10) node[below] {$(1,1,\emptyset,21)$};
\draw (4,-10) node[fill,circle,scale=0.5pt]{};
\draw (4,-10) node[below] {$(\emptyset,1,21,1)$};
\draw (6,-10) node[fill,circle,scale=0.5pt]{};
\draw (6,-10) node[below] {$(1,\emptyset,1,21)$};
\draw (8,-10) node[fill,circle,scale=0.5pt]{};
\draw (8,-10) node[below] {$(1,1,21,\emptyset)$};
\draw (3.5,-11) node {$A^{132,413265}$};
\draw (-2,-10) -- (-2,-9) -- (-0,-10) -- (5,-9) -- (6,-10) -- (1,-9) -- (-2,-10);
\draw (-2,-9) -- (4,-10) -- (8,-9);
\draw (1,-9) -- (2,-10);
\draw (8,-9) -- (8,-10);
\end{tikzpicture}}\caption{The posets of the embeddings of 132 in 413265 and the union $A^{132,413265}$ of their interiors.}\label{fig:exA} \end{figure}

Using our decomposition we build posets  from embeddings in the following~way:
\begin{defn}
Given an embedding $\eta\in E^{\sigma,\pi}$ define the poset $P(\eta):=[\hat{\eta}_1,\hat{\pi}_1]\times\cdots\times[\hat{\eta}_t,\hat{\pi}_t]$ and $$A^{\sigma,\pi}:=\bigcup_{\eta\in \widehat{E}^{\sigma,\pi}}P(\eta)^o,$$ where $P(\eta)^o$ denotes the interior of $P(\eta)$, that is, $P(\eta)$ with the top and bottom elements removed.
\end{defn}
\begin{ex}
Consider $[132,413265]$ and let $\eta_1,\eta_2,\eta_3$ and $\eta_4$ be the embeddings listed in $\widehat{E}^{\sigma,\pi}$ in Example~\ref{ex:3}. Then $\hat{\pi}=(1,1,21,21)$ and $\hat{\eta}_1=(\emptyset,1,21,\emptyset)$, $\hat{\eta}_2=(1,\emptyset,\emptyset,21)$, $\hat{\eta}_3=(\emptyset,1,\emptyset,21)$ and $\hat{\eta}_4=(\emptyset,\emptyset,1,21)$. See Figure \ref{fig:exA} for $P(\eta_i)$ and $A^{132,413265}$.
\end{ex}

The poset $A^{\sigma,\pi}$ consists of the elements $\hat{\eta}$ for all $\eta\in\widehat{E}^{\lambda,\pi}$ and all $\lambda\in(\sigma,\pi)$. Therefore, we define a surjective poset map $f$ from $A^{\sigma,\pi}$ to $(\sigma,\pi)$ in the following~way:
\begin{defn}\label{defn:map}
Let $f:A^{\sigma,\pi}\rightarrow(\sigma,\pi)$ be the map which maps all elements $\hat{\eta}$, where $\eta\in\widehat{E}^{\lambda,\pi}$, to $\lambda$. 
\end{defn}

\begin{ex}
If $[132,413265]$ and $\hat{\eta}=(1,\emptyset,1,21)$ then $\eta=400265\in \widehat{E}^{2143,\pi}$, so $f(\hat{\eta})=2143$. 
\end{ex}

\subsection{The Topology of a Poset}
We study the topology of a poset by constructing a simplicial complex from the poset in the following way:
\begin{defn}
Let $P$ be a poset. A \emph{chain} in $P$ is a totally ordered subset $\{z_1<\cdots<z_t\}$. The \emph{order complex} of $P$, denoted $\Delta(P)$, is the simplicial complex whose vertices are the elements of $P$ and whose faces are the chains of~$P$.
\end{defn}

When we refer to the order complex of an interval $[\sigma,\pi]$ we mean the order complex of the interior~$(\sigma,\pi)$, which we denote $\Delta(\sigma,\pi)$. 

\begin{figure}[h]\centering
\begin{tikzpicture}[scale=0.65]
\draw (-2,0) node[fill,circle,scale=0.5pt]{};
\draw (-2,0) node[left] {4123};
\draw (0,0) node[fill,circle,scale=0.5pt]{};
\draw (0,0) node[left] {2341};
\draw (2,0) node[fill,circle,scale=0.5pt]{};
\draw (2,0) node[right] {1234};
\draw (-2,2) node[fill,circle,scale=0.5pt]{};
\draw (-2,2) node[left] {45123};
\draw (0,2) node[fill,circle,scale=0.5pt]{};
\draw (0,2) node[left] {34512};
\draw (2,2) node[fill,circle,scale=0.5pt]{};
\draw (2,2) node[right] {23451};
\draw (-1,4) node[fill,circle,scale=0.5pt]{};
\draw (-1,4) node[left] {456123};
\draw (1,4) node[fill,circle,scale=0.5pt]{};
\draw (1,4) node[right] {345612};
\draw (-2,0) -- (-2,2) -- (-1,4) -- (0,2) -- (0,0) -- (2,2) -- (1,4) -- (0,2);
\draw (2,2) -- (2,0);
\draw (5,1) node[fill,circle,scale=0.5pt]{};
\draw (5,1) node[below] {4123};
\draw (5,4) node[fill,circle,scale=0.5pt]{};
\draw (5,4) node[above] {45123};
\draw (7,2.5) node[fill,circle,scale=0.5pt]{};
\draw (7,1.7) node {456123};
\draw (9,4) node[fill,circle,scale=0.5pt]{};
\draw (9,4) node[above] {34512};
\draw (9,1) node[fill,circle,scale=0.5pt]{};
\draw (8.7,1) node[below] {2341};
\draw (11,2.5) node[fill,circle,scale=0.5pt]{};
\draw (11.75,2.5) node[above] {345612};
\draw (11,0) node[fill,circle,scale=0.5pt]{};
\draw (11,0) node[below] {23451};
\draw (13,1) node[fill,circle,scale=0.5pt]{};
\draw (13,1) node[right] {1234};
\draw[pattern=north west lines,pattern color=blue!60,draw=black,thick] (7,2.5) -- (5,1) -- (5,4) -- (7,2.5) -- (9,4) -- (9,1) -- (7,2.5);
\draw[pattern=north west lines,pattern color=blue!60,draw=black,thick] (9,1) -- (9,4) -- (11,2.5) -- (9,1) -- (11,0) -- (11,2.5) -- (9,1);
\draw[pattern=north west lines,pattern color=blue!60,draw=black,thick] (11,0) -- (13,1) -- (11,2.5) -- (11,0);
\end{tikzpicture}
\caption{Left: Hasse diagram of $(123,4567123)$. Right: The order complex $\Delta(123,4567123)$.}\label{fig:ordercomplex}
\end{figure}
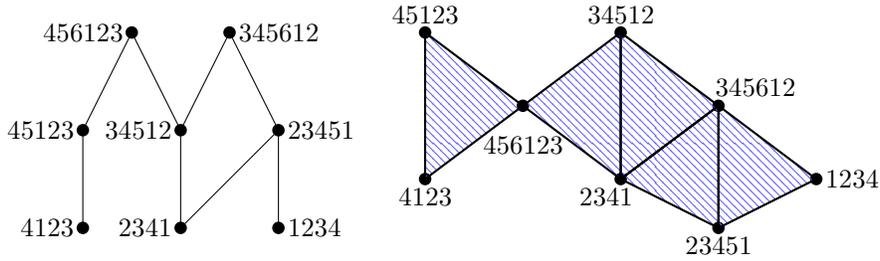

\begin{ex}
Consider the interval $I=[123,4567123]$. An example of a chain in $(123,4567123)$ is $4123<456123$. The order complex and Hasse diagram of~$I$ are given in Figure \ref{fig:ordercomplex}.
\end{ex}

We refer to a poset and its order complex interchangeably, so a topological property of a poset refers to that property of its order complex. For further background on order complexes and poset topology in general see \cite{Wac07}. 

We can use the order complex of $[\sigma,\pi]$ to calculate $\mu(\sigma,\pi)$ due to the following formula, which is an application of the Philip Hall Theorem and the Euler-Poincar\'e formula for the reduced Euler characteristic, see \cite[Section~1.2]{Wac07}:
\begin{equation}\label{eq:Euler}
\mu(\sigma,\pi)=\tilde{\chi}(\Delta(\sigma,\pi))=\sum_{i=-1}^{|\pi|-|\sigma|}(-1)^i \tilde{\beta}_i(\Delta(\sigma,\pi)),
\end{equation}
where $\tilde{\chi}$ is the reduced Euler characteristic and $\tilde{\beta}_i$ is the $i$'th reduced Betti number, that is, the rank of the $i$'th reduced homology group. Therefore, by calculating the homology of $[\sigma,\pi]$ we can compute the M\"obius function. For example, if we can show that $\Delta(\sigma,\pi)$ is contractible this implies $\mu(\sigma,\pi)=0$, and if $\Delta(\sigma,\pi)$ and $\Delta(\alpha,\beta)$ are homotopically equivalent then $\mu(\sigma,\pi)=\mu(\alpha,\beta)$.

The first explicit results on the topology of intervals of permutations appear in \cite{McSt13} and \cite{Smith13}.

\section{The Main Result}\label{sec:result}
We use the map $f$ in Definition \ref{defn:map} to calculate the M\"obius function of $[\sigma,\pi]$ by calculating $\mu(A^{\sigma,\pi})$ and the effect on the M\"obius function when applying~$f$. First we compute $\mu(A^{\sigma,\pi})$. Given a set~$A$ of posets the M\"obius function of the union of~$A$ can be calculated using the following inclusion-exclusion formula, which can be seen as a consequence of the inclusion-exclusion formula for the Euler characteristic and Equation \eqref{eq:Euler}:
\begin{equation}\label{eq:incexGen}\mu\left(\bigcup_{a\in A}a\right)=\sum_{\substack{S\subseteq A\\S\not=\emptyset}}(-1)^{|S|-1}\mu\left(\bigcap_{a\in S}a\right),\end{equation}
For more background on this see \cite{Nar74}. Applying Equation \eqref{eq:incexGen} to $A^{\sigma,\pi}$ gives:
\begin{equation}\label{eq:incex}\mu(A^{\sigma,\pi})=\sum_{\substack{S\subseteq \widehat{E}^{\sigma,\pi}\\S\not=\emptyset}}(-1)^{|S|-1}\,\mu(\bigcap_{\eta\in S}P(\eta)^o).\end{equation}
To calculate this we need to know the M\"obius function of the intersections $\cap_{\eta\in S} P(\eta)^o$. Note that when calculating the M\"obius function of a poset we add a top and bottom element. Therefore, a contractible intersection has M\"obius function 0, an empty intersection has M\"obius function $-1$ and $\mu(P(\eta)^o)=~\mu(\hat{\eta},\hat{\pi})$.

\begin{lem}\label{lem:inter}
If $S\subseteq \widehat{E}^{\sigma,\pi}$ and $|S|>1$ then: $$\mu(\bigcap_{\eta\in S}P(\eta)^o)=\begin{cases}-1,&\mbox{ if } S\in \EZ^{\sigma,\pi}\\0,&\mbox{ otherwise }\end{cases}.$$
\begin{proof}
Let $\hat{\pi}=(\hat{\pi}_1,\ldots,\hat{\pi}_t)$ and define the \emph{join} of $S$ to be $\vee S=(\max_{\eta\in S}(\hat{\eta}_1),\\\ldots,\max_{\eta\in S}(\hat{\eta}_t))$. The join is well defined because for each $i$ the set $\{\hat{\eta}_i\,|\,\eta\in S\}$ forms a chain, so there is an $\hat{\eta}_i$ that contains all others. The join of $S$ is the smallest element containing every embedding in $S$, so it is the bottom element of the intersection $I=\bigcap_{\eta\in S}P(\eta)^o$. Therefore, if $\vee S<\hat{\pi}$ then $I$ is contractible and so has M\"obius function $0$, otherwise $\vee S=\hat{\pi}$ so $I$ is empty and thus has M\"obius function $-1$. If $\vee S=\hat{\pi}$ this implies that every letter of $\pi$ is non-zero for some $\eta\in S$, that is, $S$ has empty intersection of zero sets, so~$S\in \EZ^{\sigma,\pi}$.
\end{proof}
\end{lem}

\begin{ex}
Consider our running example of $[132,413265]$. If $S^1=\{013200,\\010065\}$ then~$013200$ decomposes to $(\emptyset,1,21,\emptyset)$ and $010065$ decomposes to $\\(\emptyset,1,\emptyset,21)$, so the join is: $$\vee S^1=(\max(\emptyset,\emptyset),\max(1,1),\max(21,\emptyset),\max(\emptyset,21))=(\emptyset,1,21,21).$$ Therefore, $\vee S^1<\hat{\pi}$ so the intersection is contractible. We can check this by looking at Figure \ref{fig:exA} where we can see that the intersection $P(013200)\cap P(010065)$ is the single point $(\emptyset,1,21,21)$, which is contractible.
\end{ex}

The following sum over $\EZ^{\sigma,\pi}$ plays in important role in our results: \begin{equation*}\displaystyle\EZ(\sigma,\pi):=\sum_{S\in \EZ^{\sigma,\pi}}(-1)^{|S|}.\end{equation*}

 Now that we know the M\"obius function of the intersections we can compute~$\mu(A^{\sigma,\pi})$.
\begin{lem}\label{lem:muA}
$$\mu(A^{\sigma,\pi})=(-1)^{|\pi|-|\sigma|}\NE(\sigma,\pi)+\EZ(\sigma,\pi).$$
\begin{proof}
We can split Equation \eqref{eq:incex}  into two parts:
\begin{equation}\label{eq:incex2}\mu(A^{\sigma,\pi})=\sum_{\eta\in \widehat{E}^{\sigma,\pi}}\mu(P(\eta)^o)+\sum_{\substack{S\subseteq \widehat{E}^{\sigma,\pi}\\|S|>1}}(-1)^{|S|-1}\mu(\bigcap_{\eta\in S}P(\eta)^o).\end{equation}
By Lemma \ref{lem:inter} the second part of the right hand side of Equation \eqref{eq:incex2} \linebreak equals~$\EZ(\sigma,\pi)$. 

By the definition of $P(\eta)$, and the identity $\mu(A\times B)=\mu(A)\mu(B)$, we know $$\mu(P(\eta)^o)=\prod_{1\le i \le t} \mu(\hat{\eta}_i,\hat{\pi}_i).$$ We know that $[\hat{\eta}_i,\hat{\pi}_i]$ is always a chain, so by the definition of normality if~$\eta$ is not normal there is some $i$ such that $|\hat{\eta}_i|\le|\hat{\pi}_i|-2$, so $\mu(\hat{\eta}_i,\hat{\pi}_i)=0$, which implies  $\mu(P(\eta)^o)=0$. If $\eta$ is normal then $|\hat{\pi}_i|-|\hat{\eta}_i|=0$ or $1$, so~$\mu(\hat{\eta}_i,\hat{\pi}_i)=1$ or $-1$, for all $i$. There are $|\pi|-|\sigma|$ parts $[\hat{\eta}_i,\hat{\pi}_i]$ with $\mu(\hat{\eta}_i,\hat{\pi}_i)=-1$, one for each zero in~$\eta$, and the remaining have $\mu(\hat{\eta}_i,\hat{\pi}_i)=1$. Therefore, $\mu(P(\eta)^o)=(-1)^{|\pi|-|\sigma|}$ for each normal embedding, so the first term in the right hand side of Equation~\eqref{eq:incex2} equals~$(-1)^{|\pi|-|\sigma|}\NE(\sigma,\pi)$.
\end{proof}
\end{lem}

We now present our formula for the M\"obius function that applies to any interval of permutations:
\begin{thm}\label{thm:main}
For any permutations $\sigma$ and $\pi$:
\begin{equation*}\label{eq:main}\mu(\sigma,\pi)=(-1)^{|\pi|-|\sigma|}\NE(\sigma,\pi)+\sum_{\lambda\in[\sigma,\pi)}\mu(\sigma,\lambda) \EZ(\lambda,\pi).\end{equation*}
\begin{proof}
We take the poset $A^{\sigma,\pi}$ and for each $\lambda\in(\sigma,\pi)$ we retract $\widehat{E}^{\lambda,\pi}$ to a point we denote $\lambda$. This transforms $A^{\sigma,\pi}$ into the interval $(\sigma,\pi)$. We need to know what effect this has on the M\"obius function of~$A^{\sigma,\pi}$. 

We work our way from the bottom to the top so we can assume that all elements below the elements of $\widehat{E}^{\lambda,\pi}$ have already been retracted and all elements above have not. Define the poset $W(\lambda):=\{\tau\in A^{\sigma,\pi}\,|\,\tau\le\eta\text{ or } \tau\ge\eta\text{ for some } \eta\in\widehat{E}^{\lambda,\pi}\}$. When we retract the elements of $\widehat{E}^{\lambda,\pi}$ to $\lambda$ we retract~$W(\lambda)$ onto a contractible poset, since in that poset the element $\lambda$ is comparable to all other elements and thus represents a cone point in the corresponding order complex. This implies the change to the M\"obius function is~$-\mu(W(\lambda))$.

To compute $\mu(W(\lambda))$ we split $W(\lambda)$ into two disjoint parts \begin{align*}W(\lambda)^<&:=\{\tau\in W(\lambda)\,|\,\tau<\eta\text{ for some } \eta\in\widehat{E}^{\lambda,\pi}\},\\W(\lambda)^\ge&:=\{\tau\in W(\lambda)\,|\,\tau\ge\eta\text{ for some } \eta\in\widehat{E}^{\lambda,\pi}\}.\end{align*} The poset $W(\lambda)^<$ is isomorphic to $(\sigma,\lambda)$ because all points below $\lambda$ have already been retracted. The poset $W(\lambda)^\ge$ is equal to $\bigcup_{\eta\in\widehat{E}^{\lambda,\pi}}(P(\eta)\setminus\hat{\pi})$ whose atoms are $\widehat{E}^{\lambda,\pi}$, so by the Crosscut Theorem, see Proposition \ref{prop:crosscut}, the M\"obius function is given by~$\mu(W(\lambda)^\ge)=\EZ(\lambda,\pi)$ (this also follows by inclusion-exclusion formula and Lemma \ref{lem:inter}).

Because every element of $\widehat{E}^{\lambda,\pi}$ lies above every element of $(\sigma,\lambda)$ this \linebreak implies $W(\lambda)=W(\lambda)^<\star W(\lambda)^\ge$, where $\star$ denotes the topological join,\linebreak so  $\mu(W(\lambda))=-\mu(W(\lambda)^<)\star\mu( W(\lambda)^\ge)$ by \cite[Theorem 10.23(2)]{Koz08}. Therefore, $$-\mu(W(\lambda))=\mu(\sigma,\lambda)\EZ(\lambda,\pi).$$ So we start with $\mu(A^{\sigma,\pi})$, given by Lemma \ref{lem:muA}, and then subtract $\mu(W(\lambda))$ for each $\lambda\in (\sigma,\pi)$, which gives the desired formula.
\end{proof}
\end{thm}

\begin{rem}\label{rem:tests}
Computer tests\footnote{The programs used are available at \url{https://github.com/JasonPSmith/perm} and the data can be provided upon request by emailing the author at jason.p.smith@strath.ac.uk.} indicate that $~95\%$ of intervals $[\sigma,\pi]$, where \linebreak$|\pi|<9$, satisfy $\mu(\sigma,\pi)=(-1)^{|\pi|-|\sigma|}\NE(\sigma,\pi)$. Thus, for these intervals the latter term in Equation \eqref{eq:main} is zero.
\end{rem}

\begin{rem}
The complexity of counting the number of normal embeddings is polynomial so in the cases where we can show that the latter term of Equation~\eqref{eq:main} equals zero we have a polynomial time formula for the M\"obius function. This is a dramatic improvement over the original recursive formula that has exponential complexity. However, computing the latter term of Equation \eqref{eq:main} also has exponential complexity.
\end{rem}

Tests show that using Equation \eqref{eq:main} is often much quicker than computing the M\"obius function using the recursive formula. When computing the M\"obius function of the rank $15$ interval $$[54123, 9\, 7\, 10\, 4\, 8\, 1\, 2\, 6\, 5\, 3\, 19 \,17\, 20\, 14 \,18\, 11\, 12 \,16\, 15\, 13],$$ the formula in Equation \eqref{eq:main} took 1.75 minutes and the recursive formula took 13.5 hours. Note that this interval has M\"obius function $-3$ but no normal embeddings so the latter term of Equation \eqref{eq:main} is nonzero in this case. Furthermore, using Equation~\eqref{eq:main} we were able to compute the M\"obius function of a rank 16 interval in~1 hour and a rank 17 interval in 6 hours. However, if $\sigma$ has a large number of occurrences in $\pi$ then using Equation \eqref{eq:main} can be quite slow. For example, if $\sigma=2413$ and $\pi= 2\, 4\, 6\, 8\, 10\, 1\, 3\, 5\, 7\, 9$ then there are 35 occurrences of~$\sigma$ in $\pi$ and $\mu(\sigma,\pi)$ can be computed in 0.06 seconds using the recursive formula but takes 15.5 hours using Equation \eqref{eq:main}.

\subsection{Poset Fibration}\label{sec:fibration}

In this subsection we present a generalisation of the argument used to prove Theorem \ref{thm:main}. We can view the pair $((\sigma,\pi),\{\widehat{E}^{\lambda,\pi}\}_{\lambda\in(\sigma,\pi)})$ as a poset fibration which makes $f$ the projection map and $A^{\sigma,\pi}$ the total space. In \cite{Bjo05} various theorems are presented which relate two posets $P$ and $Q$ linked by a poset fibration~$f$ satisfying a certain condition, see Theorem 2.5 of \cite{Bjo05} for the most general form of this condition. However, our poset fibration does not always satisfy this condition, for example the condition is not true on the interval~$[1,456123]$. We present a result with a different condition on the poset fibration that generalises the argument in the proof of Theorem \ref{thm:main}, this result can be seen as an application of \cite[Corollary 3.2]{Wal81}. We let $f^*$ and $f^{-1}$ denote the  image and preimage of $f$, respectively.

\begin{prop}\label{prop:mob2.5}
Let $f:P\rightarrow Q$ be a surjective poset map such that $f^*(P_{<p})=Q_{<q}$, for any $q\in Q$ and $p\in f^{-1}(q)$. Then
$$\mu(Q)=\mu(P)+\sum_{q\in Q}\mu(Q_{<q})\mu(f^{-1}(Q_{\ge q})).$$
\begin{proof}
We begin with $P$ and for each $q\in Q$ we retract $f^{-1}(q)$ to a single point and observe the effect this retraction has on the M\"obius function of $P$. We do this inductively from the bottom to the top, so when considering $q\in Q$ we assume all points in $P_{<p}$, for all $p\in f^{-1}(q)$, have been retracted.

To calculate the effect the retraction has on the M\"obius function of the poset we consider
$$W(q):=\{p\in P\,|\, p<\lambda\text{ or }p\ge\lambda\text{ for some }\lambda\in f^{-1}(q)\}.$$ When we retract $f^{-1}(q)$ to a point we retract $W(q)$ to a contractible poset, which implies the change to the M\"obius function is $-\mu(W(q))$. We can rewrite~$W(q)$ in the following way:
\begin{align*}W(q)&=\bigcup_{p\in f^{-1}(q)}P_{<p}\star P_{\ge p}=\bigcup_{p\in f^{-1}(q)}Q_{<q}\star P_{\ge p}\\&=Q_{<q}\star\bigcup_{p\in f^{-1}(q)} P_{\ge p}=Q_{<q}\star f^{-1}(Q_{\ge q}).\end{align*}
We can replace $P_{<p}$ with $Q_{<q}$ because our induction assumption is that $P_{<p}$ has been retracted and our condition of the proposition is $f^*(P_{<p})=Q_{<q}$. Therefore, $-\mu(W(q))=\mu(Q_{<q})\mu(f^{-1}(Q_{\ge q}))$ and summing over all $q\in Q$ completes the proof.
\end{proof}
\end{prop}

\begin{rem}
An interesting question is whether Proposition \ref{prop:mob2.5} can be generalised to show homotopy equivalence. Also, is there a more general condition that encompasses the conditions of both Proposition~\ref{prop:mob2.5} and Theorem~2.5 of \cite{Bjo05}?
\end{rem}

\section{Applications}\label{sec:app}
By Theorem \ref{thm:main} we know the M\"obius function is linked to the number of normal embeddings, which depend on the adjacencies.
\begin{lem}\label{lem:adj}
The average total number of letters in the tails of adjacencies in a permutation of length~$n$ is $2(\frac{n-1}{n})$.  In particular, when $n$ tends to infinity the average number of letters in the tails of adjacencies tends to 2.
\begin{proof}
 Note first that $k=1$ cannot be in the tail of an increasing
  adjacency and $n$ cannot be in the tail of a decreasing adjacency.
 For $k>1$ the number of permutations of length $n$ in
  which $k$ is in the tail of an increasing adjacency is $(n-1)!$,
  because these are precisely all permutations of the letters
  $1,2,\ldots,n$ where $(k-1)k$ is regarded as a single letter.  So
  the probability that a letter $k>1$ is in the tail of an increasing
  adjacency is $(n-1)!/n!=1/n$. Likewise, the probability that a
  letter $k<n$ is in the tail of a decreasing adjacency is
  $1/n$. Therefore, the probability that $k$ is in the tail of an
  adjacency is
\begin{equation}\label{eq:tail}\begin{cases}\frac{1}{n},&\mbox{ if } k=1 \mbox{ or }  n\\\frac{2}{n},&\mbox{ otherwise}\end{cases}.\nonumber\end{equation}
Summing over all letters $k=1,\ldots,n$ completes the proof.
\end{proof}
\end{lem}

An embedding in a permutation $\pi$ is likely to be normal if there is only a small proportion of letters in the tails of the adjacencies of $\pi$. Therefore, Lemma~\ref{lem:adj} indicates that the proportion of embeddings of $\lambda$ in a random permutation $\pi$ that are normal increases as the length of $\pi$ increases. If a permutation has no adjacencies of size $\ell>1$ then all embeddings will be normal, the proportion of such permutations tends to~$1/e^2$ as the length of the permutations increase, see~\cite{OEISA002464}.


By Remark \ref{rem:tests} we suspect that the second part of Equation \eqref{eq:main} vanishes for a significant proportion of intervals. A consequence of Proposition 3.3 in~\cite{Smith14} is that if $\sigma$ and $\pi$ have the same number of descents then the second part of Equation \eqref{eq:main} vanishes. Note that although the definition of normal embeddings in \cite{Smith14} does not consider decreasing adjacencies, it is equivalent to Definition~\ref{defn:normal} when the number of descents is fixed. To see this note that if $\sigma$ and $\pi$ have the same number of descents then any letters in $\pi$ that form a decreasing adjacency, and thus a descent, must be nonzero in all embeddings otherwise that descent would not be in $\sigma$, contradicting the assumption that $\sigma$ and $\pi$ have the same number of descents.

One route to simplifying Equation \eqref{eq:main} is answering the following question:
\begin{que}\label{que:sum}
Given an interval $[\sigma,\pi]$, for which $\lambda\in[\sigma,\pi)$ is $\EZ(\lambda,\pi)$ nonzero?
\end{que}

One case where $\EZ(\lambda,\pi)=0$ is when $\EZ^{\lambda,\pi}=\emptyset$, which leads us to the following definition and proposition. Our notation here follows from the idea of interval blocks in \cite{SteTen10}. 
\begin{defn}
We say an interval $[\sigma,\pi]$ has a \emph{single block} if there exists some~$i$ such that $\eta_i=0$ for any $\eta\in \widehat{E}^{\sigma,\pi}$. That is, there is a letter in $\pi$ that is not contained in any of the occurrences in $\widehat{E}^{\sigma,\pi}$. 

We say an interval is \emph{single} if for all $\lambda\in[\sigma,\pi)$ the interval $[\lambda,\pi]$ has a single block.
\end{defn}


\begin{prop}\label{prop:single}
If $[\sigma,\pi]$ is single then $\mu(\sigma,\pi)=(-1)^{|\pi|-|\sigma|}\NE(\sigma,\pi)$.
\begin{proof}
If $[\lambda,\pi]$ has a single block then $\EZ^{\lambda,\pi}$ must be empty which implies $\EZ(\lambda,\pi)=~0$. Therefore, if $[\sigma,\pi]$ is single then $\EZ(\lambda,\pi)=0$ for all $\lambda\in[\sigma,\pi)$, combining this with Equation~\eqref{eq:main} completes the proof.
\end{proof}
\end{prop}

Intervals that contain a disconnected subinterval of rank at least 3 are non-shellable, as shown by Bj\"orner in \cite{Bjo80}, and thus not amenable to some of the elegant methods of topological combinatorics, see \cite{McSt13} for further  background. In the rest of the paper we consider a particular type of interval $[\alpha,\beta]$ that is known to be disconnected and show that $\EZ(\alpha,\beta)$ is nonzero for these intervals. Whether there is a topological ``reason" for $\EZ(\alpha,\beta)$  being nonzero in these cases we don't know. 

We consider decomposable permutations and write them in the form \linebreak$\pi_1\oplus\cdots\oplus\pi_n$ where each~$\pi_i$, which we call a \emph{component} of $\pi$, is indecomposable. Consider a permutation $\pi=\pi_1\oplus\cdots\oplus\pi_\ell$, where $\pi_i=\pi_1$ for any $i\in \{1,\ldots,\ell\}$, and $\lambda\le\pi$ obtained from $\pi$ by removing $k$ of the components from this sequence, where $\ell>k\ge1$. The interval $[\lambda,\pi]$ is disconnected, which follows from Lemma 4.2 in~\cite{McSt13}. Intervals of this form are the cause of the extra conditions in the formulas for the M\"obius function that appear in~\cite{BJJS11} and~\cite{McSt13}.

\begin{defn}
Given an indecomposable permutation $\lambda$ let $$\lambda^n:=\underbrace{\lambda\oplus\cdots\oplus\lambda}_{\times n}.$$
\end{defn}

Before we continue let us introduce a very useful result for computing the M\"obius function of lattices known as the Crosscut Theorem. We denote the join of $A\subseteq X$ as $\vee A:=\min\{x\in X\,|\,x\ge a\text{ for all } a\in A\}$:

\begin{prop}\label{prop:crosscut}(Crosscut Theorem, see \cite[Corollary~3.9.4]{Sta97})
Let $L$ be a finite lattice with top element $\hat{1}$ and bottom element $\hat{0}$. Let $X$ be a subset of $L$ such that $\hat{0}\not\in X$ and for all $s\in L$, $s\not=\hat{0}$, there is some $t\in X$ such that $s\ge t$. Then
$$\mu(\hat{0},\hat{1})=\sum_{\substack{A\subseteq X\\\vee A=\hat{1}}}(-1)^{|A|}.$$
\end{prop}
The Crosscut theorem is traditionally used to compute the M\"obius function of a lattice, but we can use it in reverse if we can represent our problem as a lattice for which we already know the M\"obius function. 

Consider the interval $[\lambda^m,\lambda^n]$, for some indecomposable permutation $\lambda$, the embeddings of $\lambda^m$ in~$\lambda^n$ can be considered as subsets of $[n]:=\{1,\ldots,n\}$ of size~$m$. So we can regard our problem as that of computing the M\"obius function of a sublattice of the Boolean lattice:
\begin{defn}
The \emph{Boolean lattice} $B_n$ is the poset of subsets of $[n]$ with the partial order given by inclusion.

Define the \emph{truncated Boolean lattice} $B_{n}^{\ge k}$ as the subposet of $B_n$ where all elements $a\in B_n$ such that $|a|< k$ are retracted to a single point $\hat{0}$. Similarly, define $B_{n}^{\le k}$ as the subposet of $B_n$ where all elements $a\in B_n$ such that $|a|> k$ are retracted to a single point $\hat{1}$.
\end{defn}

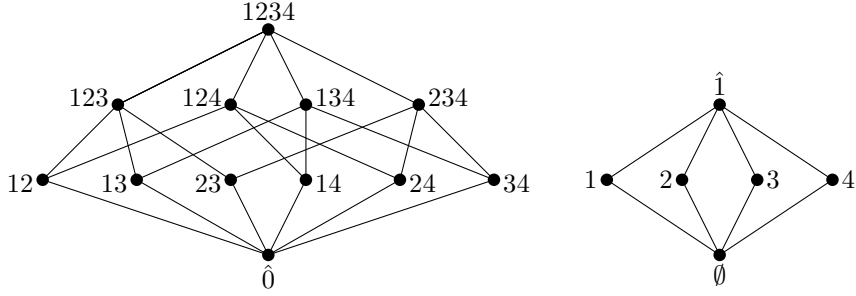
\begin{figure}[ht]\centering\begin{tikzpicture}
\draw(0,0) -- (-3,1) -- (-2,2) -- (0,3);
\draw(0,0) -- (-1.75,1) -- (-2,2) -- (0,3);
\draw(0,0) -- (0.5,1) -- (-0.5,2) -- (0,3);
\draw(0,0) -- (-0.5,1) -- (-2,2) -- (0,3);
\draw(0,0) -- (1.75,1) -- (2,2) -- (0,3);
\draw(0,0) -- (3,1) -- (0.5,2) -- (0,3);
\draw(3,1) -- (2,2);
\draw(1.75,1) -- (-0.5,2);
\draw(-0.5,1) -- (2,2);
\draw(0.5,1) -- (0.5,2);
\draw(-1.75,1) -- (0.5,2);
\draw(-3,1) -- (-0.5,2);
\draw (0,3) node[above] {$1234$};
\draw (0,3) node[fill,circle,scale=0.5pt]{};
\draw (-2,2.1) node[left] {$123$};
\draw (-2,2) node[fill,circle,scale=0.5pt]{};
\draw (-0.5,2.1) node[left] {$124$};
\draw (-0.5,2) node[fill,circle,scale=0.5pt]{};
\draw (0.5,2.1) node[right] {$134$};
\draw (0.5,2) node[fill,circle,scale=0.5pt]{};
\draw (2,2.1) node[right] {$234$};
\draw (2,2) node[fill,circle,scale=0.5pt]{};
\draw (-3,.95) node[left] {$12$};
\draw (-3,1) node[fill,circle,scale=0.5pt]{};
\draw (-1.75,.95) node[left] {$13$};
\draw (-1.75,1) node[fill,circle,scale=0.5pt]{};
\draw (-0.5,.95) node[left] {$23$};
\draw (-0.5,1) node[fill,circle,scale=0.5pt]{};
\draw (0.5,.95) node[right] {$14$};
\draw (0.5,1) node[fill,circle,scale=0.5pt]{};
\draw (1.75,.95) node[right] {$24$};
\draw (1.75,1) node[fill,circle,scale=0.5pt]{};
\draw (3,.95) node[right] {$34$};
\draw (3,1) node[fill,circle,scale=0.5pt]{};
\draw (0,0) node[below] {$\hat{0}$};
\draw (0,0) node[fill,circle,scale=0.5pt]{};
\draw(6,0) -- (4.5,1) -- (6,2);
\draw(6,0) -- (5.5,1) -- (6,2);
\draw(6,0) -- (6.5,1) -- (6,2);
\draw(6,0) -- (7.5,1) -- (6,2);
\draw (6,2) node[above] {$\hat{1}$};
\draw (6,2) node[fill,circle,scale=0.5pt]{};
\draw (4.5,1) node[left] {$1$};
\draw (4.5,1) node[fill,circle,scale=0.5pt]{};
\draw (5.5,1) node[left] {$2$};
\draw (5.5,1) node[fill,circle,scale=0.5pt]{};
\draw (6.5,1) node[right] {$3$};
\draw (6.5,1) node[fill,circle,scale=0.5pt]{};
\draw (7.5,1) node[right] {$4$};
\draw (7.5,1) node[fill,circle,scale=0.5pt]{};
\draw (6,0) node[below] {$\emptyset$};
\draw (6,0) node[fill,circle,scale=0.5pt]{};
\end{tikzpicture}
\caption{The truncated Boolean lattices $B_{4}^{\ge2}$ and $B_{4}^{\le1}$.}\label{fig:partialboolean}
\end{figure}
We take the notation for a truncated Boolean lattice from \cite[Section~3.2.1]{Wac07}. See Figure~\ref{fig:partialboolean} for examples of truncated Boolean lattices. The embeddings of~$\lambda^m$ in $\lambda^n$ can be viewed as the atoms of the truncated Boolean lattice~$B_n^{\ge m}$, so using the Crosscut theorem we can compute $\EZ(\lambda^m,\lambda^n)$ by computing~$\mu(B_n^{\ge m})$. The M\"obius function of a Boolean lattice is given by~$\mu(B_n)=(-1)^n$, see \cite[Section~3]{Rota64}. We can use this to compute the M\"obius function of the truncated Boolean lattice.

\begin{lem}\label{lem:parBool}
The M\"obius function of a truncated Boolean lattice is given by: $$\mu(B_{n}^{\le k})=(-1)^{k-1}\dbinom{n-1}{k} \hskip30pt \text{ and }\hskip30pt\mu(B_n^{\ge k})=(-1)^{n-k-1}\dbinom{n-1}{k-1}.$$
\begin{proof}
First consider $B_{n}^{\le k}$. For each element $\lambda\in B_{n}^{\le k}$, with $|\lambda|=\ell$, the interval~$[\emptyset,\lambda]$ is isomorphic to the boolean lattice $B_\ell$, therefore $\mu(\emptyset,\lambda)=(-1)^\ell$. There are $\binom{n}{\ell}$ elements in $B_{n}^{\le k}$ with size $\ell$, for $0\le\ell\le k$. To compute $\mu(B_{n}^{\le k})$ we need to sum all elements and negate, we do this by summing over $\ell$. We can then apply an identity on the alternating sum of binomial coefficients, a proof of which can be found in Section 0 of \cite{Kle63}, this gives:
$$\mu(B_{n}^{\le k})=-\sum_{\ell=0}^{k}(-1)^{\ell}\binom{n}{\ell}=(-1)^{k-1}\binom{n-1}{k}.$$
Note that the lattice $B_n^{\ge k}$ is isomorphic to $(B_{n}^{\le n-k})^*$, the dual of $B_{n}^{\le n-k}$, that is, the lattice with the partial order reversed. Therefore, $\mu(B_n^{\ge k})=\mu((B_{n}^{\le n-k})^*)=\mu(B_{n}^{\le n-k})$ which completes the proof.
\end{proof}
\end{lem}

We can now present our result for the interval $[\lambda^m,\lambda^n]$:
\begin{prop}\label{prop:countEmb}
Let $\lambda$ be an indecomposable permutation, of length $\ell>1$, and consider the interval $[\lambda^m,\lambda^n]$. Then:
$$\EZ(\lambda^m,\lambda^n)=(-1)^{n-m-1}\dbinom{n-1}{m-1}.$$
\begin{proof}
We can consider each embedding of $\lambda^m$ in $\lambda^n$ as a subset $a\subseteq\{1,\ldots,n\}$ with $|a|=m$. Therefore, the embeddings correspond to the atoms of the lattice~$B_n^{\ge m}$. So we can apply the Crosscut theorem and Lemma~\ref{lem:parBool} to complete the proof.
\end{proof}
\end{prop}

Finally, we present two results showing that, in certain cases, we can prepend or append components to both permutations $\sigma$ and $\pi$ without changing the value of $\EZ(\sigma,\pi)$. This allows us to begin with an interval for which we know the value of $\EZ(\sigma,\pi)$ and build up to larger intervals.

\begin{prop}\label{prop:prepend2}
Consider a pair of indecomposable permutations $\alpha$ and $\lambda$ such that~$[\alpha,\lambda]$ is single. Then:
$$\EZ(\alpha\oplus\lambda^m,\alpha\oplus\lambda^n)=\EZ(\lambda^m,\lambda^n)=\EZ(\lambda^m\oplus\alpha,\lambda^n\oplus\alpha).$$
\begin{proof}
First we show that $\EZ(\alpha\oplus\lambda^m,\alpha\oplus\lambda^n)=\EZ(\lambda^m,\lambda^n)$. Let $\sigma=\lambda^m$, $\pi=\lambda^n$,  $\bar{\sigma}=\alpha\oplus\sigma$ and $\bar{\pi}=\alpha\oplus\pi$. Furthermore, let $\pi_i$ be the $i$'th copy of $\lambda$ in~$\pi$. Similarly define $\sigma_i$, $\bar{\sigma}_i$ and  $\bar{\pi}_i$ and let~$\bar{\pi}_0=\bar{\sigma}_0=\alpha$.

Given an embedding $\phi$ of $\sigma$ in $\pi$, let $r(\phi)$ be the index of the first component of $\pi$ that $\phi$ embeds in, minus $1$. Each embedding $\phi$ of $\sigma$ in $\pi$ can be extended to an embedding of $\bar{\sigma}$ in $\bar{\pi}$ by choosing where to embed $\bar{\lambda}_0$ in $\bar{\pi}$. We can embed $\bar{\lambda}_0$ in~$\bar{\pi}_0$ or in each $\pi_i$ in $w:=|\hat{E}^{\alpha,\lambda}|$ different ways, for all $0<i\le r(\phi)$. Therefore, there are $1+r(\phi)w$ ways to extend $\phi$.

We can extend a set $S\in \EZ^{\sigma,\pi}$ to a set $\bar{S}\in\EZ^{\bar{\sigma},\bar{\pi}}$ by extending each embedding $\phi\in S$. In fact we can extend each $\phi$ multiple times to create~$\bar{S}$ by adding new copies of $\phi$ extended in different ways.
Note that because $S\in \EZ^{\sigma,\pi}$ every letter is nonzero for at least one embedding of $S$. So there exists at least one embedding $\eta\in S$ with $r(\eta)=0$. Therefore, the only way to extend $\eta$ is to embed $\bar{\sigma}_0$ in $\bar{\pi}_0$, which implies that any extension of $S$ is in $\EZ^{\bar{\sigma},\bar{\pi}}$.

We claim that every set of $\EZ^{\bar{\sigma},\bar{\pi}}$ can be obtained by extending a set of $\EZ^{\sigma,\pi}$. Suppose for a contradiction that $S\in\EZ^{\bar{\sigma},\bar{\pi}}$ cannot be obtained by extending a set of $\EZ^{\sigma,\pi}$. This implies that there is a $\bar{\pi}_i$, with $i>0$, which is not embedded in by $\lambda$ for any embedding in $S$. Moreover, because $[\alpha,\lambda]$ is single, if $\lambda$ is not embedded in $\bar{\pi}_i$, then there is a letter of $\bar{\pi}_i$ that is zero for all embeddings in $S$. This implies $S\not\in\EZ^{\bar{\sigma},\bar{\pi}}$, which gives a contradiction.

So we can compute $\EZ(\bar{\sigma},\bar{\pi})$ by considering each set $S\in\EZ^{\sigma,\pi}$ and how it can be extended. If we extend each element of $S$ in exactly one way, then the cardinality of $S$ does not change. However, for each additional extended form we add we must increase the cardinality by one, which changes the parity. Therefore, for each element $\phi\in S$, there are $1+r(\phi)w$ ways to extend $\phi$ and we can choose~$k$ of these for any $k=1,\ldots,1+r(\phi)w$ and this adds $k-1$ elements to the set. This gives us the following formula:

\begin{equation}\label{eq:prepend}
\EZ(\bar{\sigma},\bar{\pi})=\sum_{S\in \EZ^{\sigma,\pi}}(-1)^{|S|}\prod_{\phi\in S}\sum_{k=1}^{1+r(\phi)w}(-1)^{k-1}\dbinom{1+r(\phi)w}{k}.
\end{equation}

The sum over $k$ on the right hand side of Equation \eqref{eq:prepend} equals $1$, so the result follows immediately.

The proof is analogous to show that $\EZ(\lambda^m,\lambda^n)=\EZ(\lambda^m\oplus\alpha,\lambda^n\oplus\alpha)$.
\end{proof}
\end{prop}

\begin{prop}\label{prop:prepend}
Consider a permutation  $\pi$ and a permutation $\alpha=\alpha_1\oplus\cdots\oplus\alpha_a$. If $\alpha_a\not\le\pi$, then:
$$\EZ(\alpha\oplus\lambda,\alpha\oplus\pi)=\EZ(\lambda,\pi).$$
\begin{proof}
In any embedding of $\alpha\oplus\lambda$ in $\alpha\oplus\pi$ we must embed $\alpha_a$ in $\alpha_a$. Therefore,~$\alpha$ must  embed in  $\alpha$. This implies that the embedding set of $\alpha\oplus\lambda$ in $\alpha\oplus\pi$ can be obtained by prepending $\alpha$ to each embedding of~$\lambda$ in~$\pi$.\linebreak Furthermore, these embeddings sets have the same zero sets, which\linebreak implies $\EZ(\alpha\oplus\lambda,\alpha\oplus\pi)=\EZ(\lambda,\pi)$.
\end{proof}
\end{prop}

We know that if $[\alpha,\lambda]$ is single, then $\EZ(\alpha,\lambda)=0$. Therefore, the following conjecture is a generalisation of Proposition \ref{prop:prepend2}.

\begin{conj}
Consider a pair of indecomposable permutations $\alpha$ and $\lambda$ such that $\EZ(\alpha,\lambda)=0$. Then:
$$\EZ(\alpha\oplus\lambda^m,\alpha\oplus\lambda^n)=\EZ(\lambda^m,\lambda^n).$$
\end{conj}


\section*{Acknowledgements} 
I would like to thank the anonymous referees for their extremely useful comments and corrections which greatly improved the paper.


 \newcommand{\noop}[1]{}


\end{document}